\documentclass[11pt]{article} \topmargin=-0.5cm \textheight 220mm
\usepackage{amssymb}

   \csname @addtoreset\endcsname{equation}{section}
    \csname @addtoreset\endcsname{figure}{section}
    \csname @addtoreset\endcsname{table}{section}

\newcounter{thanksnum}
\def\thanksnumber#1
{\setcounter{thanksnum}{\value{footnote}}\setcounter{footnote}{#1}%
                     \addtocounter{footnote}{-1}\footnotemark
                     \setcounter{footnote}{\value{thanksnum}}}
\def\newtheoremz#1{\@ifnextchar[{\@othmz{#1}}{\@nthmz{#1}}}

\def\@nthmz#1#2{%
\@ifnextchar[{\@xnthmz{#1}{#2}}{\@ynthmz{#1}{#2}}}

\def\@xnthmz#1#2[#3]{\expandafter\@ifdefinable\csname #1\endcsname
{\@definecounter{#1}\@addtoreset{#1}{#3}%
\expandafter\xdef\csname the#1\endcsname{\expandafter\noexpand
  \csname the#3\endcsname \@thmcountersepz \@thmcounterz{#1}}%
\global\@namedef{#1}{\@thmz{#1}{#2}}\global\@namedef{end#1}{\@endtheoremz}}}

\def\@ynthmz#1#2{\expandafter\@ifdefinable\csname #1\endcsname
{\@definecounter{#1}%
\expandafter\xdef\csname the#1\endcsname{\@thmcounterz{#1}}%
\global\@namedef{#1}{\@thm{#1}{#2}}\global\@namedef{end#1}{\@endtheoremz}}}

\def\@othmz#1[#2]#3{\expandafter\@ifdefinable\csname #1\endcsname
  {\global\@namedef{the#1}{\@nameuse{the#2}}%
\global\@namedef{#1}{\@thmz{#2}{#3}}%
\global\@namedef{end#1}{\@endtheoremz}}}

\def\@thmz#1#2{\refstepcounter
    {#1}\@ifnextchar[{\@ythmz{#1}{#2}}{\@xthmz{#1}{#2}}}

\def\@xthmz#1#2{\@begintheoremz{#2}{\csname the#1\endcsname}\ignorespaces}
\def\@ythmz#1#2[#3]{\@opargbegintheoremz{#2}{\csname
       the#1\endcsname}{#3}\ignorespaces}

\def\@thmcounterz#1{\noexpand\arabic{#1}}
\def\@thmcountersepz{.}
\def\@begintheoremz#1#2{ \trivlist \item[\hskip \labelsep{\bf #1\ #2}]}
\def\@opargbegintheoremz#1#2#3{ \trivlist
      \item[\hskip \labelsep{\bf #1\ #2\ (#3)}]}
\def\@endtheoremz{\endtrivlist}

\newtheorem{theorem}{Theorem}[section]

\newtheorem{lemma}[theorem]{Lemma}

\newtheorem{remark}[theorem]{Remark}

\def\defi{\stackrel{{\scriptscriptstyle \Delta}}{=}}

\def\a{\alpha}

\def\O{\Omega}
\def\q{q}

\def\F{{\cal F}}
\def\w{\widehat}

\def\esssup{\mathop{\rm ess\, sup}}

\def\R{{\bf R}}
\def\E{{\bf E}}
\def\P{{\bf P}}

\def\b{\beta}

\def\t{\theta}
\def\oo{\bar}

\def\p{\partial}
\def\G{\Gamma}

\newcommand{\be}{\begin{equation}}
\newcommand{\ee}{\end{equation}}
\newcommand{\bd}{\begin{displaymath}}
\newcommand{\ed}{\end{displaymath}}
\newcommand{\ba}{\begin{array}{ll}}
\newcommand{\ea}{\end{array}}
\newcommand{\baa}{\begin{eqnarray}}
\newcommand{\eaa}{\end{eqnarray}}
\newcommand{\baaa}{\begin{eqnarray*}}
\newcommand{\eaaa}{\end{eqnarray*}}   \font\sm=cmr10

\def\Q{{\cal Q}}

\title{Optimal replication of random vectors by ordinary integrals}\author{ Nikolai Dokuchaev
\\ {\sm  Department of Mathematics \& Statistics, Curtin
University, GPO Box U1987, Perth,}\\ {\sm  6845 Western Australia.
Email N.Dokuchaev at curtin.edu.au}}
\begin{document}
\maketitle
\begin{abstract}
We consider a  problem of replication of random vectors by ordinary
integrals in the setting when  a underlying random variable is
generated by a Wiener process. The goal is to find an optimal
adapted process such that its cumulative integral at a fixed
terminal time matches this variable. The optimal process has to be
minimal in an integral norm.
\end{abstract}
\section{Introduction}
It is well known that
 random variables generated by a Wiener process can be represented via
 stochastic integrals, as is stated by the classical Clark--Ocone--Haussmann theorem.
 This result leaded to  the theory of
 backward stochastic differential equations and the martingale pricing
 method in Mathematical Finance.

We consider a problem of replication of random variables by ordinary
integrals.  The goal is to find an optimal adapted process such that
its cumulative integral at a fixed terminal time matches this
variable without error. The optimal process has to be minimal in an
integral norm. Explicit solution of this problem is found.

\section{The problem setting and the main result}
Consider a standard probability space $ (\Omega,\F, \P)$ and
standard $d$-dimensional Wiener process $w(t)$ (with $w(0)=0$) which
generates the filtration $\F_t={\sigma\{w(r):0\leq r\leq t\}}$
augmented by all the $\P$-null sets in $\F$. For $p\ge 1$ and $q\ge
1$, we denote by $L_{pq}^{n\times m}$ the class of random processes
$v(t)$ adapted to $\F_t$ with values in $\R^{n\times m}$ and such
that $\E\left(\|v(\cdot)\|_{L_q(0,T)}\right)^p<+\infty$. We denote
by $|\cdot|$ the Euclidean norm for vectors and Frobenius (i.e.,
Eucledean) norm for matrices.

Let $f$ be a $\F_\t$-measurable random vector,
 $f\in
L_2(\O,\F_T,\P;\R^n)$.  By the Martingale Representation Theorem,
there exists an unique $k_f\in L_{22}^{n\times d}$ such that \baaa
f=\E f+\int_0^{T}k_ f(t)dw(t).\label{Clarcfccc} \eaaa
 We assume that there exits $\t\in(0,T)$
such that \baaa \esssup_{t\in[\t,T]}\E
|k_f(t)|^2<+\infty.\label{kf}\eaaa
\par
\def\q{g}
 Let $\q:[0,T)\to\R$ be a given
measurable function such there exists $c>0$ and $\a\in (0.5,1)$ such
that \baa 0<\q(t)\le c(T-t)^{\a},\quad \q(t)^{-1}\le
c(1+(T-t)^{-\a}),\quad t\in[0,T).\label{GG} \eaa An example of such
a function is $\q(t)=1$ for $t<T-\tau$, $\q(t)=(T-t)^\a$ for $t\ge
T-\tau$, where $\tau\in(0,T)$ can be any number.

Let $U$ be the set of all processes  from $L_{21}^{n\times 1}$ such
that $\E\int_0^T\q(t)|u(t)|^2dt<+\infty$.

\par
 Let  $\G(t)$ be measurable
matrix process in $\R^{n\times n}$, such that $\G(t)=\q(t)G(t)$,
where  $G(t)>0$ is a symmetric positively defined matrix such that
the processes $G(t)$ and $G(t)^{-1}$ are both bounded. Clearly,
$\E\int_0^Tu(t)^\top\G(t)u(t)dt<+\infty$ for $u\in U$.
\par
Let $a\in\R^n$,  and let $b\in\R^{n\times n}$ be a non-degenerate
matrix.
\par
Consider the problem
 \baa&\hbox{Minimize}\quad
&\E\int_0^Tu(t)^\top\G(t) u(t)\,dt\quad\hbox{over}\quad u\in U\label{optim4}\\
&\hbox{ subject to} &\ba \frac{d x}{d t}(t) =A x(t)+ b u(t),\quad
t\in( 0,T)\\x(0)=a,\quad x(T)=f\quad \hbox{a.s.}\ea\label{sys} \eaa
 Note
that this problem is  a modification of a stochastic control problem
with terminal contingent claim. These problems were studied
intensively in the setting that involve backward stochastic
differential equations (BSDE); a first problem of this type was
introduced by Dokuchaev and Zhou (1999). In this setting,  a
non-zero diffusion coefficient is presented in the evolution
equation for the plant process as an auxiliary control process. Our
setting is different: a non-zero diffusion coefficient is not
allowed.
 Our setting is different  since the non-zero diffusion term is not allowed. Problem
(\ref{optim4})-(\ref{sys}) is a linear quadratic control problem.
However, it has a potential to be extended on control problems of a
general type,  similarly to the theory of controlled backward
stochastic differential equation.
\begin{lemma}\label{lemma1} Let $k_\mu(t)=R(t)^{-1}k_f(t)$, where
\baaa R(s)\defi \int_s^TQ(t)dt, \qquad Q(t)= e^{A(T-t)}
b\G(t)^{-1}b^\top e^{A^\top(T-t)}. \eaaa
 Then  $k_\mu(\cdot)\in L_{22}^{n\times d}$.
\end{lemma}
\begin{theorem}\label{ThM} Problem
(\ref{optim4})-(\ref{sys}) has a unique solution in $U$. This
solution is defined as  $\w u(\cdot)=\G(t)^{-1}b^\top
e^{A^\top(T-t)}\mu(t)$, where $\mu(t)$ a path-wise continuous
process such that \baaa \mu(t)=\oo\mu+\int_0^tk_\mu(s)dw(s),\eaaa
where $\oo\mu=R(0)^{-1}(\E f-q)$ and $q=e^{AT}a$.
\end{theorem}
\begin{remark}{\rm
Restrictions (\ref{GG}) on the choice of $\G(t)=\q(t)G(t)$ mean that
the penalty for the large size of $u(t)$ vanishes as $t\to T$. Thus,
we don't exclude fast growing $u(t)$ as $t\to T$ such that $u(t)$ is
not square integrable. This is why we select the class $U$ of
admissible controls to be wider than $L_{22}^{n\times 1}$. In
Dokuchaev (2010), a related result was obtained for a simpler case
when it was required to ensure that $x(T)=\E\{f|\F_\t\}$ for some
$\t<T$. In this setting, the exact match could be achieved only  for
$\F_\t$-measurable $f$; the optimal solution was found to be a
square integrable process.}
\end{remark}
\section{Proofs}
\def\|{|}
{\it Proof of Lemma \ref{lemma1}}. By the assumptions, we have that
$Q(t)= g(t)^{-1}\Q(t)$, where
 \baaa \Q(t)=e^{A(T-t)}
bG(t)^{-1}b^\top e^{A^\top(T-t)} \eaaa is a bounded matrix,
 $g(t)^{-1}\ge
c^{-1}(T-t)^{-\a}$ for $t\in[0,T)$. For a matrix $M=M^\top\ge 0$,
set \baaa \rho(M)=\inf_{x\in\R^n:\ |x|=1}x^\top Mx.\eaaa We have
that $\zeta=\inf_{s\in[0,T]}\rho(\Q(s))>0$ and \baaa \rho(R(t))\ge
\int_t^T\rho(Q(s))ds\ge c^{-1}\int_t^T(T-s)^{-\a}\rho(\Q(s))ds\ge
-\frac{\zeta}{ c}\frac{(T-s)^{1-\a}}{1-\a}\Bigl|_{s=t}^{s=T}. \eaaa
Hence \baa \rho(R(t))\le\frac{\zeta}{c} \frac{(T-t)^{1-\a}}{1-\a}.
\label{zm} \eaa It follows from (\ref{zm}) that \baaa
\|R(t)^{-1}\|\le C\frac{(1-\a)}{(T-t)^{1-\a}}, \eaaa for some
constant $C>0$ that is defined by $\zeta,c$ and $n$. Hence \baaa
\int_\t^T\left\|R(t)^{-1}\right\|^2dt\le C^2(1-\a)^2
\int_\t^T\frac{1}{(T-t)^{2-2\a}}dt= C^2(1-\a)^2
\frac{(T-\t)^{2\a-1}}{2\a-1}<+\infty. \eaaa It follows that \baaa
&&\E\int_0^T|k_\mu(t)|^2dt\le\E\int_0^T\left\|R(t)^{-1}\right\|^2|k_f(t)|^2dt\\&&=\E\int_0^\t
\left\|R(t)^{-1}\right\|^2|k_f(t)|^2dt+\E\int_\t^T\left\|R(t)^{-1}\right\|^2|k_f(t)|^2dt
\nonumber\\&&
\le\sup_{t\in[0,\t]}\left\|R(t)^{-1}\right\|^2\E\int_0^\t|k_f(t)|^2dt+\int_\t^T\left\|R(t)^{-1}\right\|^2\E|k_f(t)|^2dt\nonumber
\\&&=\sup_{t\in[0,\t]}\left\|R(t)^{-1}\right\|^2\E\int_0^\t|k_f(t)|^2dt+\esssup_{t\in[\t,T]}\E|k_f(t)|^2\int_\t^T\left\|R(t)^{-1}\right\|^2dt
\label{finall} \eaaa and \baa
\E\int_0^T|k_\mu(t)|^2dt<+\infty.\label{final} \eaa
This completes
the proof of  Lemma \ref{lemma1}. $\Box$
\par {\it Proof of Theorem
\ref{ThM}}. By the definition of $U$, it follows that $x(T)\in
L_2(\O,\F_T,\P;\R^n )$ for any $u\in U$.
 Let the function $L(u,\mu):U\times L_2(\O,\F_T,\P;\R^n)\to\R$ be defined as
\baaa L(u,\mu)\defi\frac{1}{2}\E\int_0^Tu(t)^\top
\G(t)u(t)\,dt+\E\mu(f-x(T)). \eaaa For a given $\mu$, consider the
following problem: \baa \hbox{Minimize $L(u,\mu)$ over $u\in
U$.}\label{min}\eaa
 We solve problem (\ref{min}) using the so-called
 Stochastic Maximum Principle that gives a necessary condition of optimality
 for stochastic control problems (see, e.g., Arkin and Saksonov
  (1979), Bensoussan (1983),  Dokuchaev and Zhou (1999),  Haussmann (1986),
 Kushner (1972)). The only $u=u_\mu$
 satisfying  these necessary conditions of
 optimality is defined by  \baa &&\w
u_{\mu}(t)=\G(t)^{-1}b^\top \psi(t),\nonumber\\
&&\psi(t)=e^{A^\top(T-t)}\mu(t),\quad \mu(t)=\E\{\mu|\F_t\}.
\label{SMP}\eaa
 Clearly, the function $L(u,\mu)$ is strictly concave in $u$,
 and this minimization problem has
 an unique solution. Therefore, this $u$ is the solution of (\ref{min}).
\par
Further,  we consider the following problem: \baa \hbox{Maximize
$L(\w u_{\mu},\mu)$ over $\mu\in
L_2(\O,\F_T,\P;\R^n)$.}\label{Lmu}\eaa
\par
 Clearly, the corresponding  $x(T)$ is
$$x(T)=\int_0^T e^{A(T-t)}b\w u_{\mu}(t)dt+e^{AT}a ,$$ and
  \baaa L(\w u_\mu,\mu)
=\frac{1}{2}\E\int_0^T\w u_{\mu}(t)^\top \G(t)\w
u_{\mu}(t)\,dt-\E\mu^\top\int_0^Te^{A(T-t)}b\w
u_{\mu}(t)dt-\E\mu^\top e^{AT}a+\E\mu^\top f. \eaaa \par We have
that
 \baaa &&\frac{1}{2}\E\int_0^T\w u_{\mu}(t)^\top \G(t)\w
u_{\mu}(t)\,dt-\E\mu^\top\int_0^Te^{A(T-t)}b\w u_{\mu}(t)dt
\\&&=\frac{1}{2}\E\int_0^T(\G(t)^{-1}b^\top \psi(t))^\top
\G(t)\G(t)^{-1}b^\top
\psi(t)\,dt-\E\mu^\top\int_0^Te^{A(T-t)}b\G(t)^{-1}b^\top
\psi(t)dt
\\&&
=\frac{1}{2}\E\int_0^T\psi(t)^\top
b\G(t)^{-1}\G(t)\G(t)^{-1}b^\top \psi(t)\,dt-\E\mu^\top\int_0^T
e^{A(T-t)}b\G(t)^{-1}b^\top \psi(t)dt
\\&&
=\frac{1}{2}\E\int_0^T\psi(t)^\top b\G(t)^{-1}b^\top
\psi(t)\,dt-\E\mu^\top\int_0^Te^{A(T-t)}b\G(t)^{-1}b^\top
\psi(t)dt
\\&&
=\frac{1}{2}\E\int_0^T[e^{A^\top(T-t)}\mu(t)]^\top
b\G(t)^{-1}b^\top
e^{A^\top(T-t)}\mu(t)\,dt-\E\mu^\top\int_0^Te^{A(T-t)}b\G(t)^{-1}b^\top
e^{A^\top(T-t)}\mu(t)dt
\\&&
=\frac{1}{2}\E\int_0^T\mu(t)^\top e^{A(T-t)} b\G(t)^{-1}b^\top
e^{A^\top(T-t)}\mu(t)\,dt-\E\mu^\top\int_0^Te^{A(T-t)}b\G(t)^{-1}b^\top
e^{A^\top(T-t)}\mu(t)dt\\&& =\frac{1}{2}\E\int_0^T\mu(t)^\top
Q(t)\mu(t)\,dt-\E\mu^\top\int_0^TQ(t)\mu(t)dt=-\frac{1}{2}\E\int_0^T\mu(t)^\top
Q(t)\mu(t)dt.
  \eaaa
We used for the last equality that
$$\E\mu^\top\int_0^TQ(t)\mu(t)dt=\E\int_0^T\mu(t)^\top
Q(t)\mu(t)dt.$$
 It follows that \baaa L(\w u_\mu,\mu)=\E
\mu^\top (f-q) -\frac{1}{2}\E\int_0^{T}\mu(t)^\top Q(t) \mu(t)\,dt.
\eaaa By the Martingale Representation Theorem, there exists
$k_\mu\in L_{22}^{n\times d}$  such that \baa
\mu=\oo\mu+\int_0^Tk_\mu(t)dw(t),\label{Clarc} \eaa where
$\oo\mu\defi \E\mu$. \index{NE NADO! It follows that \baaa L(\w
u_\mu,\mu)=\E\mu^\top f -\E\mu^\top q
-\frac{1}{2}\E\int_0^{T}\mu(t)^\top Q(t) \mu(t)\,dt. \eaaa }
  We have that \baaa
L(\w u_\mu,\mu)\!\!&=&\!\!\oo\mu^\top (\oo f-q)
-\frac{1}{2}\oo\mu^\top R(0) \oo\mu\,
-\frac{1}{2}\E\int_0^Tdt\int_0^tk_\mu(\tau)^\top Q(t)
k_\mu(\tau)\,d\tau+\E\int_0^{\t}k_\mu(t)^\top k_ f(t)\,dt\\
&=&\oo\mu^\top (\oo f-q) -\frac{1}{2}\oo\mu^\top R(0) \oo\mu\,
-\frac{1}{2}\E\int_0^Tk_\mu(\tau)^\top R(\tau)
k_\mu(\tau)\,d\tau+\E\int_0^{\t}k_\mu(t)^\top k_f(t)\,dt,
 \eaaa
 Clearly, the solution of
 problem (\ref{Lmu}) is uniquely defined by  (\ref{Clarc}) with
 \baa
\oo\mu=R(0)^{-1}(\oo f-q),\quad k_\mu(t)=R(t)^{-1}k_f(t).
\label{kmu} \eaa \

\par
By Lemma \ref{lemma1}, it follows from (\ref{final}) that
$\sup_{t\in[0,T]}\E|\mu(t)|^2<+\infty$.

 Let show that $\w u_\mu\in
U$ for any $\mu$. We have that \baaa \E\left(\int_0^T|\w
u_\mu(t)|dt\right)^2\le
C_1\E\left(\int_0^T\|\G(t)^{-1}\||\mu(t)|dt\right)^2 \le
C_2\sup_{t\in[0,T]}\E|\mu(t)|^2 \int_0^Tg(t)^{-1}dt<+\infty.\eaaa In
addition, \baaa \E\int_0^Tg(t) |u_\mu(t)|^2dt\le
C_3\E\int_0^Tg(t)|\G(t)^{-1}\mu(t)|^2dt\le
C_4\E\int_0^Tg(t)^{-1}|\mu(t)|^2dt\\
\le C_4\sup_{t\in[0,T]}\E|\mu(t)|^2 \int_0^Tg(t)^{-1}dt<+\infty.
\eaaa Here $C_i>0$ are constants defined by $A,b,n$, and $T$. Hence
$\w u_\mu\in U$.
 \par
 We found that  $\sup_\mu \inf_u L(u,\mu)$ is achieved for $(\w u_\mu,\mu)$ defined by
  (\ref{Clarc}), (\ref{kmu}), (\ref{SMP}).
  We have that $L(u,\mu)$ is strictly concave in $u\in U$ and
 affine in $\mu\in L_2(\O,\F,\P,\R^n)$. In addition, $L(u,\mu)$ is continuous in $u\in L_{22}^{n\times 1}$ given $\mu\in L_2(\O,\F,\P,\R^n)$, and
 $L(u,\mu)$
 is continuous in $\mu\in L_2(\O,\F,\P,\R^n)$ given $u\in U$.
By Proposition 2.3 from  Ekland and Temam  (1999), Chapter VI, p.
175, it follows that
 \baa
 \inf_{u\in U}\sup_\mu L(u,\mu)= \sup_\mu \inf_{u\in U} L(u,\mu).
\label{infsup} \eaa Therefore, $(u_\mu,\mu)$ defined by
(\ref{Clarc}), (\ref{kmu}), (\ref{SMP}) is the unique saddle point
for (\ref{infsup}).
\par
  Let  $U_f$ be the set of all $u(\cdot)\in U$ such that (\ref{sys}) holds.  It is easy to
see that \baaa\inf_{u\in U_f} \frac{1}{2}\E\int_0^Tu(t)^\top
\G(t)u(t)\,dt =\inf_{u\in U}\sup_\mu
  L(u,\mu),\label{infsup1}\eaaa
 and any solution  $(u,\mu)$  of
 (\ref{infsup}) is such that $u\in U_f$. It follows that
 $u_\mu\in U_f$  and it is the optimal solution for problem
 (\ref{optim4})-(\ref{sys}).
  Then the proof of Theorem \ref{ThM} follows.
  \section{Example of calculation of $u$}
Consider a model where $f=F(y(T)$, where $y(t)$ satisfies Ito
equation \baaa dy(t)=a( y(t),t)   dt +\b (  y(t),t )   dw(t). \eaaa
Here $a(x,t): \R^{n}\times \R\to\R^{n}$, $\b (x,t): \R^{n} \times \R
\to \R^{n\times n}$ are measurable bounded functions  such that the
derivative $\p \b (x,t) / \p x $ is bounded, $b(x,t)= \frac{1}{2}\b
(x,t)\,\b (x,t)^\top\ge \delta I>0$ for all $x,t$, where $\delta>0$,
$I$ is the unit matrix.
\par
Theorem \ref{ThM} can be applied as the following.

Assume first that $a(x,t)\equiv 0$ then Theorem \ref{ThM} ensures
that  $x(T)=f$ and that $\E\int_0^T\G(t)u(t)^2dt$ is minimal for
\baa u(t)=\G(t)^{-1}b^\top
e^{A^\top(T-t)}\Bigl[R(0)^{-1}(H(y(0),0)-x(0))+\int_0^tR(s)^{-1}\frac{\p
H}{\p x}(y(s),s)dy(s)\Bigr].\label{u4}\eaa Here $H$ is the solution
of the Cauchy problem for parabolic equation \baa &&\frac{\p H}{\p
t}(x,t)+\sum_{i,j=1}^n b_{ij}(x,t) \frac{\p ^2 H}{\p x_i \p
x_j}(x,t) + \sum_{i=1}^n a_{i}(x,t)\frac{\p  H}{\p x_i
}(x,t)=0,\quad t<T, \nonumber\\ &&H(x,T)=F(x).\label{parab} \eaa
 Here $b_{ij}, a_i, x_i$ are the components of the
matrix $b$ and vectors $f$, $x$.
 It can be noted that
$H(x,t)=\E\left\{F(y(T))|y(t)=x\right\}$.
\par
Assume now that $a(\cdot)\neq 0$. We still have that
$\int_0^Tu(t)dt=f$
 for $u(t)$ defined by (\ref{u4})-(\ref{parab}); in this case,
 $H(x,t)=\E_Q\left\{F(y(T))|y(t)=x\right\}$,
where $\E_Q$ is the expectation is under a probability measure $Q$
such that the process $y(t)$ is a martingale under $Q$. By Girsanov
Theorem, this measure exists, it is equivalent to the original
measure $\P$ and unique under our assumptions on $a$ and $\b$. In
this case, the value $\E\int_0^T\G(t)u(t)^2dt$ is not minimal over
$u$ anymore. Instead, $\E_Q\int_0^T\G(t)u(t)^2dt$ is minimal. This
still means that deviations of $u$ are minimal but in a different
metric. It can be also noted that the definition of the class $U$
for the original measure has to be adjusted for the new measure $Q$,
with the expectations $\E$ replaced by $\E_Q$.
\par
This model could have practical applications in a number of
settings, where it is required to ensure that a controlled
differentiable process $x(t)$ matches a random vector $f$ generated
by a uncontrolled stochastic  process $y(t)$. For instance, $x(t)$
may represent a controlled path of an anti-aircraft missile, and the
process $y(t)$ may represent the driving force of an aircraft such
that $f=F(y(T)$ represents  the aircraft coordinates at time $T$.
\subsection*{Acknowledgment}  This work  was supported by  ARC grant of Australia DP120100928 to the author.

\end{document}